\documentclass[a4paper,10pt]{article}
\usepackage{latexsym,amssymb,amsmath,amsthm,amsxtra,xypic}
\usepackage{stmaryrd}
\usepackage{amsfonts}
\usepackage{bbm}
\usepackage{amscd}
\usepackage{mathrsfs}
\usepackage[dvips]{graphicx}
\usepackage[all]{xy}

\newcommand{\allowpagebreak}
\allowdisplaybreaks

\newtheorem{thm}{Theorem}[section]
\newtheorem{lem}[thm]{Lemma}
\newtheorem{cor}[thm]{Corollary}
\newtheorem{pro}[thm]{Proposition}

\theoremstyle{definition}
\newtheorem{defi}[thm]{Definition}

\newcommand {\emptycomment}[1]{}

\newcommand{\pf}{\noindent{\bf Proof.}\ }

\def\DD{D}

\newcommand{\lam}{\lambda}
\newcommand{\si}{\sigma}

\newcommand{\ot}{\otimes}

\newcommand{\gl}{\mathfrak {gl}}

\newcommand{\g}{T}
\newcommand{\h}{V}
\newcommand{\frkg}{T}
\newcommand{\frkh}{V}

\newcommand{\dlam}{_\lambda}

\newcommand{\E}{\mathrm{E}}
\renewcommand{\L}{\mathcal{L}}

\newcommand{\Hom}{\mathrm{Hom}}

\newcommand{\Ker}{\mathrm{Ker}}

\newcommand{\End}{\mathrm{End}}
\newcommand{\Ext}{\mathrm{Ext}}

\newcommand{\ad}{\mathrm{ad}}

\newcommand{\id}{\mathrm{id}}

\begin{document}

\allowdisplaybreaks

\title{Notes on Cohomologies of Lie Triple Systems}

\author{Tao Zhang}

\date{}
\maketitle

\begin{abstract}
  The cohomology theory of Lie triple systems in the sense of Yamaguti is studied by means of cohomology of Leibniz algebras in the sense of Loday.
  The notion of Nijenhuis operators  for Lie triple systems is introduced to describe trivial deformations.
  We also study the abelian extensions of Lie triple systems in details.

  \noindent{\textsl{Mathematics Subject Classification 2010:}} 17A40, 17A30, 17B56.

  \noindent{\textsl{Key Words and Phrases:}} Lie triple systems, cohomology, deformations, Nijenhuis operators, abelian extensions.
\end{abstract}

\section{Introduction}

Lie triple systems first appeared in E.Cartan's study of totally geodesic submanifolds and symmetric spaces in \cite{Ca}.
As an algebraic structure, the concept of Lie triple systems was introduced by N. Jacobson in \cite{Ja1} and studied later by W.G. Lister in \cite{Li}.

A Lie triple system consists of a vector space $\frkg$ together with a trilinear
map $[\cdot,\cdot,\cdot]: \g\times\g\times\g\to \g$ satisfying
\begin{eqnarray*}
&&[x_1, x_1,x_2]=0,\\
&&[x_1, x_2,x_3]+[x_2, x_3,x_1]+[x_3, x_1,x_2]=0,\\
&&[x_1, x_2, [y_1, y_2, y_3]] = [[ x_1, x_2, y_1], y_2, y_3] + [y_1, [ x_1, x_2, y_2], y_3] + [y_1, y_2, [ x_1, x_2, y_3]],
\end{eqnarray*}
for all $x_i, y_i \in \g$.

Some years later, the representation and cohomology theory of Lie triple systems was established by K. Yamaguti in \cite{Ya},
and studied in \cite{Ha,HP} and  \cite{KT} from a different point of view.
K. Yamaguti's work can be described as follows.

Let $\frkg$ be a Lie triple system and $V$ be a vector space. Then $(V, \theta)$ is called a representation
of $\frkg$ if the following conditions are satisfied,
\begin{itemize}
\item[$\bullet$]{\rm(R1)}\quad $[\DD(x_1,x_2),\theta(y_1,y_2)]=\theta((x_1,x_2)\circ (y_1,y_2))$,
\item[$\bullet$]{\rm(R2)}\quad $\theta(x_1,[y_1, y_2, y_3])= \theta (y_2, y_3)\theta(x_1,y_1) - \theta (y_1, y_3)\theta(x_1,y_2) + \DD (y_1, y_2)\theta(x_1,y_3)$,
\end{itemize}
where  $\theta$ is a linear map from $\ot^2\frkg$ to $\End(V)$ , $\DD(x_1,x_2):=\theta(x_2,x_1)-\theta(x_1,x_2)$
and $[\DD(x_1,x_2),\theta(y_1,y_2)]$ is the commutator of $\DD(x_1,x_2)$ and $\theta(y_1,y_2)$.

Let $(V,\theta)$ be a representation of a Lie triple system $\frkg$. Denote by $C^{2n+1}(\g,V)$
the vector space spanned by all $2n+1$-linear mapping $\omega$ of $T\times\cdots\times T$ into $V$ satisfying
$$\omega(x_1, x_2, \cdots , x_{2n-2}, x, x, y) = 0,$$
and
$$\omega(x_1, x_2,\cdots , x_{2n-2}, x, y, z) + \omega(x_1, x_2,\cdots , x_{2n-2}, y, z, x) + \omega(x_1, x_2,\cdots , x_{2n-2}, z, x, y) = 0.$$

The Yamaguti coboundary operator $\delta^{2n-1}: C^{2n-1}(T, V)\to C^{2n+1}(T, V)$ is defined as the following:
\begin{eqnarray}\label{coboudary}
\notag&&\delta^{2n-1} \omega (x_1, x_2,\cdots, x_{2n+1})\\
\notag&=&\theta(x_{2n}, x_{2n+1})\omega (x_1, x_2,\cdots, x_{2n-1})-\theta(x_{2n-1}, x_{2n+1})\omega (x_1, x_2,\cdots, x_{2n-2}, x_{2n})\\
\notag&&+\sum_{k=1}^n (-1)^{n+k}\DD(x_{2k-1}, x_{2k})\omega (x_1, x_2,\cdots,\widehat{x_{2k-1}}, \widehat{x_{2k}},\cdots, x_{2n+1})\\
&&+\sum_{k=1}^n\sum^{2n+1}_{j=2k+1}(-1)^{n+k+1}
\omega(x_1, x_2,\cdots,\widehat{x_{2k-1}}, \widehat{x_{2k}},\cdots, [x_{2k-1}, x_{2k}, x_j],\cdots, x_{2n+1}),
\end{eqnarray}
where \,$\widehat{{}}$\, denotes omission. It is proved in \cite{Ya} that $\delta^{2n+1}\circ\delta^{2n-1}=0$ and hence we get the cohomology group
$H^{\bullet}(\frkg, V)=Z^{\bullet}(\frkg, V)/B^{\bullet}(\frkg, V)$ where $Z^{\bullet}(\frkg, V)$ is the space of cocycles
and $B^{\bullet}(\frkg, V)$ is the space of coboundaries.

One may be puzzled at first at K. Yamaguti's theory of Lie triple systems since it is intrinsic.
Consider the recent development of Leibniz $n$-algebra in \cite{CLP}.
Recall that a Leibniz  algebra (see \cite{Loday}) is a vector space $\L$ together with a bracket
$\circ: \L  \times \L  \to \L $ satisfying the following left Leibniz identity
\begin{eqnarray*}
x\circ (y\circ z)=(x\circ y)\circ z+y\circ(x\circ z),
\end{eqnarray*}
for all $x,y,z\in \L$.
The authors in \cite{CLP} showed that for a Leibniz $n$-algebra, the fundamental set is a Leibniz  algebra,
and the cohomology of Leibniz $n$-algebras is defined to be the cohomology of the associated Leibniz  algebras.

We use the same idea of that paper and prove that,
given a Lie triple system $T$ and a representation $(V,\theta)$,
we have a Leibniz algebra on the fundamental set $\L=\ot^2 \frkg$ and a representation of $\L$ on the space $\Hom(\g,V)$,
and the cohomology theory of a Lie triple system defined by Yamaguti in \cite{Ya}
is in fact the cohomology theory of Leibniz algebra defined by J.-L. Loday and T. Pirashvili in \cite{Loday}.
The main result of the first part of this paper is
\begin{thm}\label{mainthm}
The Yamaguti cohomology complex of a Lie triple system $\frkg$ with coefficients in $\frkh$ coincides with the Loday-Pirashvili cohomology complex of Leibniz algebra $\L$ with coefficients in $\Hom(\frkg,\frkh)$ with the same coboundary operator.
\end{thm}

The second part of this paper is about deformations and extensions of Lie triple systems.
We study infinitesimal deformations of Lie triple systems.
The notion of Nijenhuis operators  for Lie triple system is introduced.
We also study the abelian extension of Lie triple systems in details.
Given a representation $(V,\theta)$, we prove that there is a one-to-one correspondence
between equivalence classes of abelian extensions of Lie triple systems $\frkg$ by $\h$
and the cohomology group $\mathbf{H}^3(\frkg, V)$.

The organization of  this paper is as follows.
In section \ref{sec:2}, we review some basic facts about Lie triple systems.
We show that given a Lie triple system $T$ and a representation $(V,\theta)$, we can get a Leibniz algebra $\L$ and a representation of $\L$ on the space
$\Hom(\g, V)$.  Using this fact, we prove Theorem \ref{mainthm}.
Low dimensional cohomologies are also spelled out in this section.
In section \ref{sec:3}, we study infinitesimal deformations of Lie triple systems.
The notion of Nijenhuis operators for Lie triple systems is introduced  to describe trivial deformations.
In the last section \ref{sec:4}, we also verify that Yamaguti's cohomology theory can be used to characterize abelian extensions of Lie triple systems.

Throughout this paper, all Lie triple systems are assumed to be over
an algebraically closed field $\mathbb{F}$ of characteristic different from 2 and 3.
The space of linear maps from $V$ to $W$ is denoted by $\Hom(V,W)$.


\section{Lie triple systems and cohomology}
\label{sec:2}

A Lie triple system consists of a vector space $\frkg$ together with a trilinear
map $[\cdot,\cdot,\cdot]: \g\times\g\times\g\to \g$ satisfying
\begin{eqnarray}
&&\label{eq:Lts01}[x_1, x_1,x_2]=0,\\
&&\label{eq:Lts02}[x_1, x_2,x_3]+[x_2, x_3,x_1]+[x_3, x_1,x_2]=0,\\
&&\label{eq:Lts03}[x_1, x_2, [y_1, y_2, y_3]] = [[ x_1, x_2, y_1], y_2, y_3] + [y_1, [ x_1, x_2, y_2], y_3] + [y_1, y_2, [ x_1, x_2, y_3]],
\end{eqnarray}
for all $x_i, y_i \in \g$. Note that \eqref{eq:Lts01} yields
\begin{eqnarray*}
&&\label{eq:Lts01'}[x_1, x_2,x_3]+[x_2, x_1,x_3]=0,
\end{eqnarray*}
which means the trilinear
map $[\cdot,\cdot,\cdot]$ is antisymmetric in the first two variables. We call \eqref{eq:Lts02} the Jacobi identity and  \eqref{eq:Lts03} the fundamental identity. Any Lie algebra considered with the trilinear map $[x_1, x_2,x_3]:=[[x_1, x_2],x_3]$ is a Lie triple system. In this case, \eqref{eq:Lts02} is
just the Jacobi identity of the Lie algebra.

A homomorphism between two Lie triple systems $T$ and $S$ is a map $\varphi:T\to S$ satisfying
$$\varphi[x_1, x_2,x_3]=[\varphi x_1, \varphi x_2,\varphi x_3].$$

Denote by $x=(x_1, x_2)$ and $\ad^L(x) y_i=[x_1, x_2, y_i]$, then the above equation \eqref{eq:Lts03} can be rewritten in the form
\begin{eqnarray}\label{eq:Jacobi-3'}
\ad^L(x)[y_1, y_2, y_3]=[\ad^L(x)y_1, y_2, y_3]+[y_1, \ad^L(x)y_2, y_3]+[y_1, y_2, \ad^L(x) y_3],
\end{eqnarray}
which means that $\ad^L(x)$ is a derivation of the map $[\cdot,\cdot,\cdot]$.

Denoted by $\L:=\ot{}^2\frkg$, which is called fundamental set.
The elements $x=(x_1,x_2)\in \ot^2 \frkg$ are called fundamental objects.
Define an operation on fundamental objects by
\begin{eqnarray}\label{eq:fundamental}
x\circ y=([x_1,x_2,y_1],y_2)+(y_1,[x_1,x_2,y_2]).
\end{eqnarray}
It is easy to prove that $\L$ is a Leibniz algebra, see \cite{CLP,DT} for more details.  We also have the following equality
\begin{eqnarray}
\ad^L(x)\ad^L(y) (w)-\ad^L(y)\ad^L(x) (w)=\ad^L(x\circ y) (w), 
\end{eqnarray}
for all $x,y,z\in\L, w\in \frkg$. Thus $\ad^L: \L \to \gl(\frkg)$ is a homomorphism of Leibniz algebras.

\emptycomment{
$$[\ad^L(x),\ad^L(y)]=\ad^L(x\circ y)$$
If we define $\ad^L(x_1, x_2)x_3=[x_1, x_2,x_3]$ and $\ad^R(x_1, x_2)x_3=[x_3,x_1, x_2]$, then by \eqref{eq:Lts01} and \eqref{eq:Lts02} we know that
\begin{eqnarray*}
\ad^L(x_1,x_2)(x_3)&=&[x_1, x_2,x_3]=-[x_2, x_3,x_1]-[x_3, x_1,x_2]\\
&=&[x_3, x_2,x_1]-[x_3, x_1,x_2]=\ad^R(x_2,x_1)(x_3)-\ad^R(x_1,x_2)(x_3),
\end{eqnarray*}
}

\emptycomment{
\begin{eqnarray*}
x\circ (y\circ z)=([x_1,x_2,[y_1,y_2,z_1]],z_2)+([y_1,y_2,z_1],[x_1,x_2,z_2])+([x_1,x_2,z_1],[y_1,y_2,z_2])+(z_1,[x_1,x_2,[y_1,y_2,z_2]]).
\end{eqnarray*}
\begin{eqnarray*}
y\circ (x\circ z)=([y_1,y_2,[x_1,x_2,z_1]],z_2)+([x_1,x_2,z_1],[y_1,y_2,z_2])+([y_1,y_2,z_1],[x_1,x_2,z_2])+(z_1,[y_1,y_2,[x_1,x_2,z_2]]).
\end{eqnarray*}
\begin{eqnarray*}
(x\circ y)\circ z=([[x_1,x_2,y_1],y_2,z_1],z_2)+(z_1,[[x_1,x_2,y_1],y_2,z_2])+([y_1,[x_1,x_2,y_2],z_1],z_2)+(z_1,[y_1,[x_1,x_2,y_2],z_2]).
\end{eqnarray*}
}

For a Leibniz algebra $\L$, a representation of $\L$ is a vector space $M$ together with two bilinear maps
$$[\cdot,\cdot]_L:\L\times M\to M \,\,\, \text{and}\,\,\,  [\cdot,\cdot]_R: M\times \L\to M$$
satisfying the following three axioms: $\forall x,y \in \L, m\in M$,
\begin{itemize}
\item[$\bullet$] {\rm(LLM)}\quad  $[x\circ y, m]_L=[x, [y, m]_L]_L-[y,[x, m]_L]_L$,
\item[$\bullet$] {\rm(MLL)}\quad  $[m, x\circ y]_R=[[m, x]_R, y]_R+[x, [m, y]_R]_L$,
\item[$\bullet$] {\rm(LML)}\quad  $[x,[m, y]_R]_L=[[x,m]_L, y]_R+[m, x\circ y]_R.$
\end{itemize}
By (MLL) and (LML) we also have
\begin{itemize}
\item[$\bullet$] {\rm(MMM)}\quad  $[[m, x]_R, y]_R+[[x,m]_L, y]_R=0.$
\end{itemize}
In fact, assume (LLM), any of (LML),(MLL) and (MMM) can be derived from the other two.

\begin{defi}\cite{Ya} Let $\frkg$ be a Lie triple system and $V$ be a vector space. Then $(V, \theta)$ is called a representation
of $\frkg$ (or a $\frkg$-module) if the following conditions are satisfied,
\begin{itemize}
\item[$\bullet$]{\rm(R1)}\quad $[\DD(x_1,x_2),\theta(y_1,y_2)]=\theta((x_1,x_2)\circ (y_1,y_2))$,
\item[$\bullet$]{\rm(R2)}\quad $\theta(x_1,[y_1, y_2, y_3])= \theta (y_2, y_3)\theta(x_1,y_1) - \theta (y_1, y_3)\theta(x_1,y_2) + \DD (y_1, y_2)\theta(x_1,y_3)$,
\end{itemize}
where  $\theta$ is a map from $\L=\frkg^{\ot2}$ to $\End(V)$ and $\DD(x_1,x_2):=\theta(x_2,x_1)-\theta(x_1,x_2)$.
\end{defi}

For example, given a Lie triple system $\frkg$, there is a natural adjoint representation of $\L$ on $T$.
The corresponding representation $\theta$ and $\DD$ is given by
\begin{eqnarray*}
\ad^R(x_1,x_2)(x_3)=[x_3,x_1,x_2]\quad\mbox{and}\quad \ad^L(x_1,x_2)(x_3)=[x_1,x_2,x_3].
\end{eqnarray*}

\emptycomment{
Note that by the Jacobi identity of $[\cdot,\cdot,\cdot]$ in $T$ we have
\begin{eqnarray*}
\ad^L(x_1,x_2)=\ad^R(x_2,x_1)-\ad^R(x_1,x_2).
\end{eqnarray*}

If we define
\begin{eqnarray*}
\ad^L(x_1,x_2)x_3=[x_1,x_2,x_3], \quad \ad^M(x_1,x_2)x_3=[x_1,x_3,x_2],
\end{eqnarray*}
then by the antisymmetricity and Jacobi identity of $[\cdot,\cdot,\cdot]$ in $T$ we have
\begin{eqnarray*}
\ad^L(x_1,x_2)=\ad^R(x_2,x_1)-\ad^R(x_1,x_2), \quad \ad^M(x_1,x_2)=-\ad^R(x_1,x_2).
\end{eqnarray*}
}

Given a Lie triple system $\frkg$ and a representation $(\frkh,\theta)$, we define maps
\begin{eqnarray}
[\cdot,\cdot]_L:\L \times \Hom(\frkg,\frkh)\to \Hom(\frkg,\frkh)\,\,\,
\text{and}\,\,\, [\cdot,\cdot]_R: \Hom(\frkg,\frkh)\times \L \to \Hom(\frkg,\frkh),
\end{eqnarray}
by
\begin{eqnarray}
\label{eq:leibniz01}{[(x_1,x_2),\phi]_L}(x_3)&=&\DD(x_1,x_2)\phi(x_3)-\phi([x_1,x_2,x_3])
\end{eqnarray}
and
\begin{eqnarray}
\notag[\phi,(x_1,x_2)]_R(x_3)
&=&\phi([x_1,x_2,x_3])-\DD(x_1,x_2)\phi(x_3)+\theta(x_1,x_3)\phi(x_2)\\
\label{eq:leibniz02}&&-\theta(x_2,x_3)\phi(x_1),
\end{eqnarray}
for all $\phi\in \Hom(\frkg,\frkh), x_i\in \frkg$. We will prove that there is a representation of the Leibniz algebra $\L$ on $\Hom(\frkg,\frkh)$.
Note that by \eqref{eq:leibniz01} and \eqref{eq:leibniz02} we have
\begin{eqnarray}
\label{eq:leibniz03}{([(x_1,x_2),\phi]_L+[\phi,(x_1,x_2)]_R)}(x_3)=\theta(x_1,x_3)\phi(x_2)-\theta(x_2,x_3)\phi(x_1).
\end{eqnarray}

\begin{pro}\label{pro:rep}
Let $\frkg$ be a Lie triple system. Then $\Hom(\frkg,\frkh)$ equipped with the above two maps
$[\cdot,\cdot]_L$ and $[\cdot,\cdot]_R$ is a representation of Leibniz algebra $\L$ on $\Hom(\frkg,\frkh)$
if and only if (R1), (R2) and the following (R3) are satisfied
\begin{itemize}
\item[$\bullet$]{\rm(R3)}\quad $[\DD(x_1,x_2),\DD(y_1,y_2)]=\DD((x_1,x_2)\circ (y_1,y_2))$.
\end{itemize}
\end{pro}

\pf For $x=(x_1,x_2), y=(y_1,y_2)\in\L$ and $y_3\in \frkg$, first we compute the equality
$$[x\circ y, \phi]_L(y_3)=[x, [y, \phi]_L]_L(y_3)-[y,[x, \phi]_L]_L(y_3).$$
By definition, the left hand side is equal to
\begin{eqnarray*}
[x\circ y, \phi]_L(y_3)=\DD(x\circ y)\phi(y_3)-\phi(\ad^L(x\circ y)y_3),
\end{eqnarray*}
and the right hand side is equal to
\begin{eqnarray*}
&&[x, [y, \phi]_L]_L(y_3)-[y,[x, \phi]_L]_L(y_3)\\
&=&\DD(x)[y, \phi]_L(y_3)-[y, \phi]_L(\ad^L(x)y_3)-\DD(y)[x, \phi]_L(y_3)+[x, \phi]_L(\ad^L(y)y_3)\\
&=&\DD(x)\DD(y)\phi(y_3)-\DD(x)\phi(\ad^L(y)y_3)-\DD(y)\phi(\ad^L(x)y_3)+\phi(\ad^L(y)\ad^L(x)y_3)\\
&&-\DD(y)\DD(x)\phi(y_3)+\DD(y)\phi(\ad^L(x)y_3)+\DD(x)\phi(\ad^L(y)y_3)-\phi(\ad^L(x)\ad^L(y)y_3)\\
&=&\DD(x)\DD(y)\phi(y_3)+\phi(\ad^L(y)\ad^L(x)y_3)-\DD(y)\DD(x)\phi(y_3)-\phi(\ad^L(x)\ad^L(y)y_3)\\
&=&[\DD(x),\DD(y)]\phi(y_3)-\phi([\ad^L(x),\ad^L(y)]y_3).
\end{eqnarray*}
Since $\ad^L: \L \to \End(\frkg)$ is a homomorphism of Leibniz algebras,
then (LLM) is valid for $[\cdot,\cdot]_L$ if and only if (R3) holds.

Now we compute the equality
$$[x, [\phi, y]_R]_L(y_3)=[\phi, x\circ y]_R(y_3)+[[x,\phi]_L, y]_R(y_3)$$
By definition, the left hand side is equal to
\begin{eqnarray*}
[x, [\phi, y]_R]_L(y_3)&=&\DD(x_1,x_2)[\phi, (y_1,y_2)]_R(y_3)-[\phi, (y_1,y_2)]_R([x_1,x_2,y_3])\\
&=&\DD(x_1,x_2)\{\phi([y_1,y_2,y_3])-\DD(y_1,y_2)\phi(y_3)\\
&&+\theta(y_1,y_3)\phi(y_2)-\theta(y_2,y_3)\phi(y_1)\}\\
&&-\phi([y_1,y_2,[x_1,x_2,y_3]])-\DD(y_1,y_2)\phi([x_1,x_2,y_3])\\
&&+\theta(y_1,[x_1,x_2,y_3])\phi(y_2)-\theta(y_2,[x_1,x_2,y_3])\phi(y_1),
\end{eqnarray*}
and the right hand side is
\begin{eqnarray*}
[\phi, x\circ y]_R(y_3)&=&[\phi,([x_1,x_2,y_1],y_2)]_R(y_3)+[\phi,(y_1,[x_1,x_2,y_2])]_R(y_3)\\
&=&\phi([[x_1,x_2,y_1],y_2,y_3])-\DD([x_1,x_2,y_1],y_2)\phi(y_3)\\
&&+\theta([x_1,x_2,y_1],y_3)\phi(y_2)-\theta(y_2,y_3)\phi([x_1,x_2,y_1])\\
&& + \phi([y_1,[x_1,x_2,y_2],y_3])-\DD(y_1,[x_1,x_2,y_2])\phi(y_3)\\
&&+\theta(y_1,y_3)\phi([x_1,x_2,y_2])-\theta([x_1,x_2,y_2],y_3)\phi(y_1),
\end{eqnarray*}
\begin{eqnarray*}
[[x,\phi]_L, y]_R(y_3)&=&[x,\phi]_L([y_1,y_2,y_3])-\DD(y_1,y_2)[x,\phi]_L(y_3)\\
&&+\theta(y_1,y_3)[x,\phi]_L(y_2)-\theta(y_2,y_3)[x,\phi]_L(y_1)\\
&=&\DD(x_1,x_2)\phi([y_1,y_2,y_3])-\phi([x_1,x_2,[y_1,y_2,y_3]])\\
&&-\DD(y_1,y_2)\{\DD(x_1,x_2)\phi(y_3)-\phi([x_1,x_2,y_3])\}\\
&&+\theta(y_1,y_3)\{\DD(x_1,x_2)\phi(y_2)-\phi([x_1,x_2,y_2])\}\\
&&-\theta(y_2,y_3)\{\DD(x_1,x_2)\phi(y_1)-\phi([x_1,x_2,y_1])\},
\end{eqnarray*}
Therefore (LML) is valid for $[\cdot,\cdot]_L$ if and only if (R1) and (R3) hold.

At last, we compute the equality
$$[[\phi, x]_R+[x,\phi]_L, y]_R(y_3)=0.$$
By \eqref{eq:leibniz01} and \eqref{eq:leibniz02} we have
$${[(x_1,x_2),\phi]_L(w)+[\phi,(x_1,x_2)]_R}(w)=\theta(x_1,w)\phi(x_2)-\theta(x_2,w)\phi(x_1),$$
thus
\begin{eqnarray*}
&&[[\phi, x]_R+[x,\phi]_L, y]_R(y_3)\\
&=&([\phi, x]_R+[x,\phi]_L)([y_1,y_2,y_3])-\DD(y_1,y_2)([\phi, x]_R+[x,\phi]_L)(y_3)\\
&&+\theta(y_1,y_3)([\phi, x]_R+[x,\phi]_L)(y_2)-\theta(y_2,y_3)([\phi, x]_R+[x,\phi]_L)(y_1)\\
&=&\theta(x_1,[y_1,y_2,y_3])\phi(x_2)-\theta(x_2,[y_1,y_2,y_3])\phi(x_1)\\
&&-\DD(y_1,y_2)\theta(x_1,y_3)\phi(x_2)+\DD(y_1,y_2)\theta(x_2,y_3)\phi(x_1)\\
&&+\theta(y_1,y_3)\theta(x_1,y_2)\phi(x_2)-\theta(y_1,y_3)\theta(x_2,y_2)\phi(x_1)\\
&&-\theta(y_2,y_3)\theta(x_1,y_1)\phi(x_2)+\theta(y_2,y_3)\theta(x_2,y_1)\phi(x_1).
\end{eqnarray*}
Therefore (MMM) is valid for $[\cdot,\cdot]_L$ and $[\cdot,\cdot]_R$ if and only if (R2) holds.
\qed

The relationship between (R1) and (R3) is as follows.
Assume (R1), then we have
\emptycomment{
\begin{itemize}
\item[$\bullet$]{\rm(R1)}\quad
$[\DD(x_1,x_2),\theta(y_1,y_2)]=\theta((x_1,x_2)\circ (y_1,y_2)))=\theta([x_1,x_2,y_1],y_2)+\theta(y_2,[x_1,x_2,y_1])$.
\end{itemize}
\begin{itemize}
\item[$\bullet$]{\rm(R1')}\quad
$[\DD(x_1,x_2),\theta(y_2,y_1)]=\theta((x_1,x_2)\circ (y_2,y_1))=\theta([x_1,x_2,y_2],y_1)+\theta(y_1,[x_1,x_2,y_2])$.
\end{itemize}
Now it is easy to see
}
\begin{eqnarray*}
[\DD(x_1,x_2),\DD(y_1,y_2)]
&=&[\DD(x_1,x_2),\theta(y_2,y_1)-\theta(y_1,y_2)]\\
&=&\theta([x_1,x_2,y_2],y_1)+\theta(y_1,[x_1,x_2,y_2])\\
&&-\theta([x_1,x_2,y_1],y_2)-\theta(y_2,[x_1,x_2,y_1])\\
&=&\DD([x_1,x_2,y_1],y_2)+\DD(y_1, [x_1,x_2,y_2])\\
&=&\DD((x_1,x_2)\circ (y_1,y_2)).
\end{eqnarray*}
Thus from (R1) we get (R3). In general,  from (R3) we can't get (R1) since $\theta$ is not an antisymmetric map.
Nevertheless, we have
\begin{cor}\label{pro:rep}
Let $\frkg$ be a Lie triple system. Then $\Hom(\frkg,\frkh)$ equipped with the above two maps
$[\cdot,\cdot]_L$ and $[\cdot,\cdot]_R$ is a representation of Leibniz algebra $\L$ on $\Hom(\frkg,\frkh)$.
\end{cor}

Now we can prove the main result of this section.

{\bf The proof of Theorem \ref{mainthm}:}
The fact that the Yamaguti cochain complex is equal to the cochain complex of the Leibniz algebra $\L$ with coefficients in $\Hom(\frkg,\frkh)$
is clear by
$$C^{2n+1}(\frkg,\frkh)=\Hom\left(\ot{}^{2n+1}\frkg,\frkh\right)
\cong\Hom\left(\ot^{n}\L,\Hom(\frkg,\frkh)\right).$$
For the Leibniz algebra $\L$ with representation on $\Hom(\frkg,\frkh)$, the coboundary operator is defined to be
\begin{eqnarray*}
&&d_{n-1}\omega(x^1,x^2,\cdots,x^n,w)\\
&=&d_{n-1}\omega(x^1,x^2,\cdots,x^n)(w)\\
&=&\sum_{k=1}^{n-1}(-1)^{k+1}[x^k,\omega(x^1,\cdots,\hat{x^k},\cdots,x^n)]_L(w)+(-1)^n[\omega(x^1,\cdots,x^{n-1}),x^n]_R(w)\\
&&+\sum_{1\leq k<l\leq n}(-1)^{k}\omega(x^1,\cdots,\hat{x^k},\cdots,x^{l-1},x^k\circ x^l,x^{k+1},\cdots, x^n)(w),
\end{eqnarray*}
for all $x^i\in \L=\ot{}^2\frkg,\ w\in \frkg$. For more details of cohomology of  Leibniz algebras, see \cite{Loday}.

Put $x^k=(x_{2k-1}, x_{2k})$, $w=x_{2n+1}$ and $[\cdot,\cdot]_L$, $[\cdot,\cdot]_R$ as in \eqref{eq:leibniz01} and \eqref{eq:leibniz02},
then we get a coboundary operator $\delta^{2n-1}=(-1)^{n+1}d_{n-1}: C^{2n-1}(T, V)\to C^{2n+1}(T, V)$ as follows:
\begin{eqnarray*}
&&\delta^{2n-1}\omega (x_1, x_2,\cdots, x_{2n+1}):=(-1)^{n+1}d_{n-1}\omega(x^1,x^2,\cdots,x^n,w)\\
&=&(-1)^{n+1}\left\{\sum_{k=1}^{n-1} (-1)^{k+1}[(x_{2k-1}, x_{2k}),\omega (x_1, x_2,\cdots,\widehat{x_{2k-1}}, \widehat{x_{2k}},\cdots, x_{2n})]_L(x_{2n+1})\right.\\
&&+(-1)^n[\omega(x_1,x_2,\cdots,x_{2n-2}),(x_{2n-1}, x_{2n})]_R(x_{2n+1})\\
&&+\left.\sum_{k=1}^n\sum^{2n}_{j=2k+1}(-1)^{k}\omega(x_1,\cdots,\widehat{x_{2k-1}}, \widehat{x_{2k}},\cdots,[x_{2k-1}, x_{2k}, x_j],\cdots, x_{2n})(x_{2n+1})\right\}\\
&=&(-1)^{n+1}\left\{\sum_{k=1}^{n-1} (-1)^{k+1}\{\DD(x_{2k-1}, x_{2k})\omega (x_1, x_2,\cdots,\widehat{x_{2k-1}}, \widehat{x_{2k}},\cdots, x_{2n+1})\right.\\
&&-\omega(x_1, x_2,\cdots,\widehat{x_{2k-1}}, \widehat{x_{2k}},\cdots, [x_{2k-1}, x_{2k}, x_{2n+1}])\}\\
&&+(-1)^n \{\omega(x_1, x_2,\cdots,[x_{2n-1}, x_{2n}, x_{2n+1}])-\DD(x_{2n-1}, x_{2n})\omega (x_1, x_2,\cdots, x_{2n+1})\\
&&+\theta(x_{2n-1}, x_{2n+1})\omega (x_1, x_2,\cdots, x_{2n-2}, x_{2n})-\theta(x_{2n}, x_{2n+1})\omega (x_1, x_2,\cdots, x_{2n-1})\}\\
&&+\left.\sum_{k=1}^n\sum^{2n}_{j=2k+1}(-1)^{k}
\omega(x_1, x_2,\cdots,\widehat{x_{2k-1}}, \widehat{x_{2k}},\cdots, [x_{2k-1}, x_{2k}, x_j],\cdots, x_{2n+1})\right\}\\
&=&\theta(x_{2n}, x_{2n+1})\omega (x_1, x_2,\cdots, x_{2n-1})-\theta(x_{2n-1}, x_{2n+1})\omega (x_1, x_2,\cdots, x_{2n-2}, x_{2n})\\
&&+\sum_{k=1}^n (-1)^{n+k}\DD(x_{2k-1}, x_{2k})\omega (x_1, x_2,\cdots,\widehat{x_{2k-1}}, \widehat{x_{2k}},\cdots, x_{2n+1})\\
&&+\sum_{k=1}^n\sum^{2n+1}_{j=2k+1}(-1)^{n+k+1}
\omega(x_1, x_2,\cdots,\widehat{x_{2k-1}}, \widehat{x_{2k}},\cdots, [x_{2k-1}, x_{2k}, x_j],\cdots, x_{2n+1}).
\end{eqnarray*}
This is exactly the Yamaguti coboundary $\delta$ as in \eqref{coboudary}.

\qed


That is why Yamaguti can define the cohomology of Lie triple systems in 1960:
the fundamental set is a Leibniz algebra, the space $\Hom(\g,V)$ is its representation,
and the Yamaguti cohomology is nothing but the cohomology of Leibniz algebras defined by J.-L. Loday and T. Pirashvili \cite{Loday} thirty-three years later!

Now we list the low dimensional coboundary operators which will be used in the following sections.
We use the coboundary operator $d$ instead of Yamaguti coboundary operator $\delta$ to remind the readers how the maps defined in
\eqref{eq:leibniz01} and \eqref{eq:leibniz02} work.
According to the above definition, a 1-cochain is a map $\nu\in\Hom(\frkg,\frkh)$,
a 3-cochain is a map
$\omega\in \Hom\left(\ot{}^2\frkg,\Hom(\frkg,\frkh)\right)=\Hom\left(\ot{}^3\frkg,\frkh\right)$,
and the coboundary operator is given by
\begin{eqnarray}
\label{eq:cobound01}d_{0}\nu(x^1,w)&=&d_{0}\nu(x^1)(w)=-[\nu,x^1]_R(w),\\
\label{eq:cobound02}d_{1}\omega(x^1,x^2,w)&=&[x^1,\omega(x^2)]_L(w)+[\omega(x^1),x^2]_R(w)-\omega(x^1\circ x^2)(w),\\
\label{eq:cobound03} \notag d_{2}\omega(x^1,x^2,x^3,w)&=&[x^1,\omega(x^2,x^3)]_L(w)-[x^2,\omega(x^1,x^3)]_L(w)-[\omega(x^1,x^2),x^3]_R(w)\\
 &&-\omega(x^1\circ x^2,x^3)(w)+\omega(x^1,x^2\circ x^3)(w)-\omega(x^2,x^1\circ x^3)(w).
\end{eqnarray}

\emptycomment{
\begin{eqnarray}
\label{eq:cobound04}  \notag
 d_{2}\omega(x^1,x^2,x^3,x^4,w)&=&[x^1,\omega(x^2,x^3,x^4)]_L(w)-[x^2,\omega(x^1,x^2,x^3)]_L(w)+[x^3,\omega(x^1,x^2,x^4)]_L(w)\\
 \notag&&+[\omega(x^1,x^2,x^3),x^4]_R(w)-\omega(x^1\circ x^2,x^3,x^4)(w)-\omega(x^2,x^1\circ x^3,x^4)(w)\\
 \notag&&-\omega(x^2,x^3,x^1\circ x^4)(w)+\omega(x^1,x^2\circ x^3,x^4)(w)+\omega(x^1,x^2\circ x^4,x^3)(w)\\
       &&+\omega(x^1,x^2,x^3\circ x^4)(w).
\end{eqnarray}
}

Put $x^1=(x_1,x_2)\in \L,\ w=x_3\in\frkg$ in the equality \eqref{eq:cobound01}, then by \eqref{eq:leibniz02} we have
\begin{eqnarray}\label{eq:1coc}
\nonumber&&d_0\nu(x_1, x_2,x_3)=\DD(x_1, x_2)\nu(x_3)-\theta(x_1,x_3)\nu(x_2)+\theta(x_2,x_3)\nu(x_1)-\nu([x_1, x_2, x_3]).
\end{eqnarray}
\begin{defi}
Let $\frkg$ be a Lie triple system and $(V, \theta)$ be a $\frkg$-module. Then a map $\nu\in\Hom(\frkg,\frkh)$
is called 1-cocycle if
\begin{eqnarray}\label{eq:0coc}
\DD(x_1, x_2)\nu(x_3)-\theta(x_1,x_3)\nu(x_2)+\theta(x_2,x_3)\nu(x_1)-\nu([x_1, x_2, x_3])=0,
\end{eqnarray}
and a map $\omega: \ot^3\frkg\to V$ is called a 3-coboudary if there exists a map $\nu\in\Hom(\frkg,\frkh)$ such that $\omega=d_0\nu$.
\end{defi}
Put $x^1=(x_1,x_2)\in \L, x^2=(y_1,y_2)\in \L$, $w=y_3\in\frkg$ in the equality \eqref{eq:cobound02}, then we have
\begin{eqnarray*}
d_{1}\omega(x_1,x_2,y_1,y_2,y_3)&=&[x_1,x_2,\omega(y_1,y_2)]_L(y_3)+[\omega(x_1,x_2),y_1,y_2]_R(y_3)\\
&&-\omega((x_1,x_2)\circ (y_1,y_2))(y_3)\\
&=&\DD(x_1,x_2)\omega(y_1,y_2,y_3)-\omega(y_1,y_2,[x_1,x_2,y_3])\\
&&+\omega(x_1,x_2,[y_1,y_2,y_3])-\DD(y_1,y_2)\omega(x_1,x_2,y_3)\\
&&+\theta(y_1,y_3)\omega(x_1,x_2,y_2)-\theta(y_2,y_3)\omega(x_1,x_2,y_1)\\
&&-\omega([x_1,x_2,y_1],y_2,y_3)-\omega(y_1,[x_1,x_2,y_2],y_3).
\end{eqnarray*}

\emptycomment{
where
\begin{eqnarray*}
[x_1,x_2,\omega(y_1,y_2)]_L(y_3)&=&\DD(x_1,x_2)\omega(y_1,y_2)(y_3)-\omega(y_1,y_2)([x_1,x_2,y_3])\\
&=&\DD(x_1,x_2)\omega(y_1,y_2,y_3)-\omega(y_1,y_2,[x_1,x_2,y_3]),
\end{eqnarray*}
\begin{eqnarray*}
[\omega(x_1,x_2),y_1,y_2]_R(y_3)&=&\omega(x_1,x_2)([y_1,y_2,y_3])-\DD(y_1,y_2)\omega(x_1,x_2)(y_3)\\
&&+\theta(y_1,y_3)\omega(x_1,x_2)(y_2)-\theta(y_2,y_3)\omega(x_1,x_2)(y_1)\\
&=&\omega(x_1,x_2,[y_1,y_2,y_3])-\DD(y_1,y_2)\omega(x_1,x_2,y_3)\\
&&+\theta(y_1,y_3)\omega(x_1,x_2,y_2)-\theta(y_2,y_3)\omega(x_1,x_2,y_1),
\end{eqnarray*}
\begin{eqnarray*}
\omega((x_1,x_2)\circ (y_1,y_2))(y_3)&=&\omega(([x_1,x_2,y_1],y_2)+(y_1,[x_1,x_2,y_2]))(y_3)\\
&=&\omega([x_1,x_2,y_1],y_2,y_3)+\omega(y_1,[x_1,x_2,y_2],y_3).
\end{eqnarray*}
}

\begin{defi}\label{def:1coc}
Let $\frkg$ be a Lie triple system and $(V, \theta)$ be a $\frkg$-module. Then a map
$\omega: \ot^3\frkg\to V$ is called 3-cocycle if $\forall x_1, x_2,x_3, y_1, y_2, y_3\in \frkg$,
\begin{eqnarray}\label{eq:1coc}
&&\omega(x_1, x_1, x_2)=0,\\
&&\omega(x_1, x_2, x_3)+ \omega(x_2, x_3, x_1)+ \omega(x_3, x_1, x_2)=0,\\
\notag&& \omega( x_1, x_2,[y_1, y_2, y_3])+\DD(x_1, x_2)\omega(y_1, y_2, y_3)\\
\notag&&=\omega([x_1, x_2, y_1], y_2, y_3) + \omega(y_1,[x_1, x_2,y_2], y_3) + \omega(y_1, y_2, [x_1, x_2,y_3])\\
\label{eq:1coc03}&&\qquad+\theta(y_2, y_3)\omega(x_1,x_2,y_1)- \theta(y_1, y_3)\omega(x_1,x_2,y_2) + \DD(y_1, y_2)\omega(x_1,x_2,y_3).
\end{eqnarray}
\end{defi}

\emptycomment{
\begin{eqnarray*}
{[(x_1,x_2),\phi]_L}(x_3)&=&\DD(x_1,x_2)\phi(x_3)-\phi([x_1,x_2,x_3]),\\
{[\phi,(x_1,x_2)]_R}(x_3)&=&\phi([x_1,x_2,x_3])-\DD(x_1,x_2)\phi(x_3)+\theta(x_1,x_3)\phi(x_2)-\theta(x_2,x_3)\phi(x_1)
\end{eqnarray*}
\begin{eqnarray*}
d_{2}\omega(x^1,x^2,x^3,w)&=&[x^1,\omega(x^2,x^3)]_L(w)-[x^2,\omega(x^1,x^3)]_L(w)-[\omega(x^1,x^2),x^3]_R(w)\\
 &&-\omega(x^1\circ x^2,x^3)(w)+\omega(x^1,x^2\circ x^3)(w)-\omega(x^2,x^1\circ x^3)(w).
\end{eqnarray*}
}

Put $x^1=(x_1,x_2)\in \L, x^2=(y_1,y_2), x^3=(z_1,z_2)\in \L$, $w=z_3\in\frkg$ in the equality \eqref{eq:cobound03}. Then we have
\begin{eqnarray*}
 &&d_{2}\omega(x_1,x_2,y_1,y_2,z_1,z_2,z_3)\\
 &=&[x_1,x_2,\omega(y_1,y_2,z_1,z_2)]_L(z_3)-[y_1,y_2,\omega(x_1,x_2,z_1,z_2)]_L(z_3)\\
 &&-[\omega(x_1,x_2,y_1,y_2),z_1,z_2]_R(z_3)-\omega((x_1,x_2)\circ (y_1,y_2),z_1,z_2)(z_3)\\
 &&+\omega(x_1,x_2,(y_1,y_2)\circ (z_1,z_2))(z_3)-\omega(y_1,y_2,(x_1,x_2)\circ (z_1,z_2))(z_3)\\
 &=&\DD(x_1,x_2)\omega(y_1,y_2,z_1,z_2,z_3)-\omega(y_1,y_2,z_1,z_2,[x_1,x_2,z_3])\\
 &&-\DD(y_1,y_2)\omega(x_1,x_2,z_1,z_2,z_3)+\omega(x_1,x_2,z_1,z_2,[y_1,y_2,z_3])\\
 &&-\omega(x_1,x_2,y_1,y_2,[z_1,z_2,z_3])+\DD(z_1,z_2)\omega(x_1,x_2,y_1,y_2,z_3)\\
 &&-\theta(z_1,z_3)\omega(x_1,x_2,y_1,y_2,z_2)+\theta(z_2,z_3)\omega(x_1,x_2,y_1,y_2,z_1)\\
 &&-\omega([x_1,x_2,y_1],y_2,z_1,z_2,z_3)-\omega(y_1,[x_1,x_2,y_2],z_1,z_2,z_3)\\
 &&+\omega(x_1,x_2,[y_1,y_2,z_1],z_2,z_3)+\omega(x_1,x_2,z_1,[y_1,y_2,z_2],z_3)\\
 &&-\omega(y_1,y_2,[x_1,x_2,z_1],z_2,z_3)-\omega(y_1,y_2,z_1,[x_1,x_2,z_2],z_3).
\end{eqnarray*}

\begin{defi}\label{def:2coc}
Let $\frkg$ be a Lie triple system and $(V, \theta)$ be a $\frkg$-module. Then a map
$\omega: \ot^5\frkg\to V$ is called 5-cocycle if $\forall x_1, x_2,y_1, y_2, z_1,z_2, z_3\in \frkg$,
\begin{eqnarray}\label{eq:1coc}
&&\omega(x_1, x_2, y_1, y_1, y_2)=0,\\
&&\omega(x_1, x_2, z_1, z_2, z_3)+ \omega(x_1, x_2, z_2, z_3, z_1)+ \omega(x_1, x_2, z_3, z_1, z_2)=0,\\
 \notag&&\DD(x_1,x_2)\omega(y_1,y_2,z_1,z_2,z_3)-\DD(y_1,y_2)\omega(x_1,x_2,z_1,z_2,z_3)\\
 \notag&&+\DD(z_1,z_2)\omega(x_1,x_2,y_1,y_2,z_3)-\theta(z_1,z_3)\omega(x_1,x_2,y_1,y_2,z_2)\\
 \notag&&+\theta(z_2,z_3)\omega(x_1,x_2,y_1,y_2,z_1)-\omega([x_1,x_2,y_1],y_2,z_1,z_2,z_3)\\
 \notag&&-\omega(y_1,[x_1,x_2,y_2],z_1,z_2,z_3)-\omega(y_1,y_2,[x_1,x_2,z_1],z_2,z_3)\\
 \notag&&-\omega(y_1,y_2,z_1,[x_1,x_2,z_2],z_3)-\omega(y_1,y_2,z_1,z_2,[x_1,x_2,z_3])\\
 \notag&&+\omega(x_1,x_2,[y_1,y_2,z_1],z_2,z_3)+\omega(x_1,x_2,z_1,[y_1,y_2,z_2],z_3)\\
       &&+\omega(x_1,x_2,z_1,z_2,[y_1,y_2,z_3])-\omega(x_1,x_2,y_1,y_2,[z_1,z_2,z_3])=0.
\end{eqnarray}
\end{defi}

\emptycomment{
Put $x^1=(x_1,x_2)\in \L, x^2=(x_3,x_4), x^3=(x_5,x_6)\in \L$, $x_7\in\frkg$ in the equality \eqref{eq:cobound03}, then we have
\begin{eqnarray*}
&&d_{1}\omega(x_1,x_2,x_3,x_4,x_5,x_6,x_7)\\
&=&[x_1,x_2,\omega(x_3,x_4,x_5,x_6)]_L(x_7)-[x_3,x_4,\omega(x_1,x_2,x_5,x_6)]_L(x_7)-[\omega(x_1,x_2,x_3,x_4),x_5,x_6]_R(x_7)\\
 &&-\omega((x_1,x_2)\circ (x_3,x_4),x_5,x_6)(x_7)+\omega(x_1,x_2,(x_3,x_4)\circ (x_5,x_6)(x_7)-\omega(x_3,x_4,(x_1,x_2)\circ x_5,x_6))(x_7)\\
 &=&\DD(x_1,x_2)\omega(x_3,x_4,x_5,x_6,x_7)-\omega(x_3,x_4,x_5,x_6,[x_1,x_2,x_7])\\
 &&-\DD(x_3,x_4)\omega(x_1,x_2,x_5,x_6,x_7)+\omega(x_1,x_2,x_5,x_6,[x_3,x_4,x_7])\\
 &&-\omega(x_1,x_2,x_3,x_4,[x_5,x_6,x_7])+\DD(x_5,x_6)\omega(x_1,x_2,x_3,x_4,x_7)\\
 &&-\theta(x_5,x_7)\omega(x_1,x_2,x_3,x_4,x_6)+\theta(x_6,x_7)\omega(x_1,x_2,x_3,x_4,x_5)\\
 &&-\omega([x_1,x_2,x_3],x_4,x_5,x_6,x_7)-\omega(x_3,[x_1,x_2,x_4],x_5,x_6,x_7)\\
 &&+\omega(x_1,x_2,[x_3,x_4,x_5],x_6,x_7)+\omega(x_1,x_2,x_5,[x_3,x_4,x_6],x_7)\\
 &&-\omega(x_3,x_4,[x_1,x_2,x_5],x_6,x_7)-\omega(x_3,x_4,x_5,[x_1,x_2,x_6],x_7).
\end{eqnarray*}
}

\section{Infinitesimal Deformations of Lie triple systems}
\label{sec:3}

In this section, we study infinitesimal deformations of Lie triple systems. We introduce the notion of Nijenhuis operators for Lie triple systems,
which is analogous to the case of ordinary Lie algebras in \cite{D,KM}.
This kind of operators give trivial deformations. For the general deformations of Lie triple systems, see \cite{KT}.

Let $\frkg$ be a Lie triple system and $\omega:\g\times\g\times\g\to\g$ be a trilinear map. Consider a $\lambda$-parametrized family of linear operations:
\begin{eqnarray*}
[x_1, x_2, x_3]_\lam&\triangleq& [x_1, x_2, x_3]+ \lambda\omega(x_1, x_2, x_3),
 \end{eqnarray*}
 where $\lambda$ is a formal variable. 

If $[\cdot,\cdot,\cdot]_\lam$ endow $\frkg$ with Lie triple system structure which is denoted by $\frkg_\lam$, then we say that $\omega$ generates a
$\lambda$-parameter infinitesimal deformation of the Lie triple system $\frkg$.

\begin{thm}\label{thm:deformation}
$\omega$ generates a $\lambda$-parameter infinitesimal deformation of the Lie triple system $\frkg$ is equivalent to
(i) $\omega$ itself defines a Lie triple system structure on $\g$ and
(ii) $\omega$ is a 3-cocycle of $\frkg$ with the coefficients in the adjoint representation.
\end{thm}

\pf
From the equality
\begin{eqnarray*}
0=[x_1, x_1, x_3]_\lam= [x_1, x_1, x_2]+ \lambda\omega(x_1, x_1, x_2),
 \end{eqnarray*}
we have
\begin{eqnarray}
\label{eq:dm01}\omega(x_1, x_1, x_2)=0.
 \end{eqnarray}

From the equality
\begin{eqnarray*}
0&=&[x_1, x_2,x_3]\dlam+[x_2, x_3,x_1]\dlam+[x_3, x_1,x_2]\dlam\\
&=&[x_1, x_2,x_3]+[x_2, x_3,x_1]+[x_3, x_1,x_2]\\
&&+\lam\{\omega(x_1, x_2,x_3)+\omega(x_2, x_3,x_1)+\omega(x_3, x_1,x_2)\},
\end{eqnarray*}
we have
\begin{eqnarray}
\label{eq:dm02}\omega(x_1, x_2,x_3)+\omega(x_2, x_3,x_1)+\omega(x_3, x_1,x_2)=0.
\end{eqnarray}

For the equality
\begin{eqnarray*}
&&[x_1, x_2, [y_1, y_2, y_3]\dlam]\dlam \\
&=&[[ x_1, x_2, y_1]\dlam , y_2, y_3]\dlam  + [y_1, [ x_1, x_2, y_2]\dlam , y_3]\dlam  + [y_1, y_2, [ x_1, x_2, y_3]\dlam]\dlam,
\end{eqnarray*}
the left hand side is equal to
\begin{eqnarray*}
&&[x_1, x_2, [y_1, y_2, y_3]+\lam\omega(y_1, y_2, y_3)]\dlam\\
&=&[x_1, x_2, [y_1, y_2, y_3]]+\lam\omega(x_1, x_2, [y_1, y_2, y_3])\\
&&+[x_1, x_2,\lam\omega(y_1, y_2, y_3)]+\lam\omega(x_1, x_2,\lam\omega(y_1, y_2, y_3))\\
&=&[x_1, x_2, [y_1, y_2, y_3]]+\lam\{\omega(x_1, x_2, [y_1, y_2, y_3])+[x_1, x_2,\omega(y_1, y_2, y_3)]\}\\
&&+\lam^2\omega(x_1, x_2,\omega(y_1, y_2, y_3)),
\end{eqnarray*}
and the right hand side is equal to
\begin{eqnarray*}
&&[[x_1, x_2, y_1]+\lam\omega(x_1, x_2, y_1), y_2, y_3]\dlam  + [y_1, [ x_1, x_2, y_2]+\lam\omega(x_1, x_2, y_2), y_3]\dlam \\
&&+ [y_1, y_2, [ x_1, x_2, y_3]+\lam\omega(x_1, x_2, y_3)]\dlam\\
&=&[[ x_1, x_2, y_1] , y_2, y_3]  + [y_1, [ x_1, x_2, y_2] , y_3] + [y_1, y_2, [ x_1, x_2, y_3]]\\
&&+\lam\{\omega([x_1, x_2, y_1],y_2, y_3) +[\omega(x_1, x_2, y_1), y_2, y_3]\\
&&\qquad+\omega(y_1,[ x_1, x_2, y_2], y_3)+[y_1,\omega(x_1, x_2, y_2), y_3]\\
&&\qquad+\omega(y_1, y_2, [ x_1, x_2, y_3])+[y_1,y_2,\omega(x_1, x_2, y_3)] \}\\
&&+\lam^2\{\omega(\omega( x_1, x_2, y_1), y_2, y_3) + \omega(y_1, \omega( x_1, x_2, y_2), y_3) + \omega(y_1, y_2, \omega( x_1, x_2, y_3))\}.
\end{eqnarray*}
Thus we have
\begin{eqnarray}
\nonumber &&\omega(x_1, x_2, [y_1, y_2, y_3])+\ad^L(x_1, x_2)\omega(y_1, y_2, y_3)\\
\nonumber &=&\omega([x_1, x_2, y_1],y_2, y_3)+\omega(y_1,[ x_1, x_2, y_2], y_3) +\omega(y_1, y_2, [ x_1, x_2, y_3])\\
\label{eq:2-coc01} &&+\ad^R(y_2, y_3)\omega(x_1, x_2, y_1)-\ad^R(y_1, y_3)\omega(x_1, x_2, y_2)+\ad^L(y_1,y_2)\omega(x_1, x_2, y_3),
\end{eqnarray}
and
\begin{eqnarray}
\nonumber &&\omega(x_1, x_2,\omega(y_1, y_2, y_3))\\
\label{eq:dm03}&=&\omega(\omega( x_1, x_2, y_1), y_2, y_3) + \omega(y_1, \omega( x_1, x_2, y_2), y_3)
+ \omega(y_1, y_2, \omega( x_1, x_2, y_3)).
\end{eqnarray}
Therefore by \eqref{eq:dm01}, \eqref{eq:dm02} and \eqref{eq:dm03} $\omega$ defines a Lie triple system structure on $\g$.
Furthermore, by \eqref{eq:2-coc01} $\omega$ is a 3-cocycle of $\frkg$ with the coefficients in the adjoint representation.
\qed


A deformation is said to be trivial if there exists a linear map $N:\frkg\to \frkg$
such that for $\varphi_\lam = \id + \lambda N$: $\frkg_\lam \to \frkg$ we have
\begin{eqnarray}
\varphi_\lam [x_1,x_2,x_3]\dlam=[\varphi_\lam x_1,\varphi_\lam x_2,\varphi_\lam x_3].
\end{eqnarray}
By definition we have
\begin{eqnarray*}
\varphi_\lam [x_1,x_2,x_3]\dlam&=&[x_1,x_2,x_3]+\lam\omega(x_1,x_2,x_3)+\lam N([x_1,x_2,x_3]+\lam\omega(x_1,x_2,x_3))\\
&=&[x_1,x_2,x_3]+\lam(\omega(x_1,x_2,x_3)+ N[x_1,x_2,x_3])+\lam^2N\omega(x_1,x_2,x_3),
\end{eqnarray*}
and
\begin{eqnarray*}
[\varphi_\lam x_1,\varphi_\lam x_2,\varphi_\lam x_3]&=&[x+\lam Nx_1,x_2+\lam Nx_2,x_3+\lam Nx_3]\\
&=&[x_1,x_2,x_3]+\lam([Nx_1,x_2,x_3]+[x_1,Nx_2,x_3]+[x_1,x_2,Nx_3])\\
&&+\lam^2([Nx_1,Nx_2,x_3]+[Nx_1,x_2,Nx_3]+[x_1,Nx_2,Nx_3])\\
&&+\lam^3[Nx_1,Nx_2,Nx_3].
\end{eqnarray*}
Thus we have
\begin{eqnarray}
\label{eq:Nijenhuis1}\omega(x_1,x_2,x_3)&=&[Nx_1,x_2,x_3]+[x_1,Nx_2,x_3]+[x_1,x_2,Nx_3]-N[x_1,x_2,x_3],\\
\label{eq:Nijenhuis2}N\omega(x_1,x_2,x_3)&=&[Nx_1,Nx_2,x_3]+[Nx_1,x_2,Nx_3]+[x_1,Nx_2,Nx_3],\\
\label{eq:Nijenhuis3}0&=&[Nx_1,Nx_2,Nx_3].
\end{eqnarray}

From the cohomology theory discussed in section \ref{sec:2}, \eqref{eq:Nijenhuis1} can be represented in terms of 1-coboundary
as $\omega=d_0N$. Moreover, it follows from \eqref{eq:Nijenhuis1} and \eqref{eq:Nijenhuis2} that $N$ must satisfy the following condition
\begin{eqnarray}\label{eq:Nijenhuis4}
\nonumber N^2[x_1,x_2,x_3]&=&N[Nx_1,x_2,x_3]+N[x_1,Nx_2,x_3]+N[x_1,x_2,Nx_3]\\
&&-([Nx_1,Nx_2,x_3]+[Nx_1,x_2,Nx_3]+[x_1,Nx_2,Nx_3]).
\end{eqnarray}
In the following, we denote by $\omega(x_1,x_2,x_3)=[x_1,x_2,x_3]_N$, then \eqref{eq:Nijenhuis2} is equivalent to
\begin{eqnarray}\label{eq:Nijenhuis4'}
N[x_1,x_2,x_3]_N&=&[Nx_1,Nx_2,x_3]+[Nx_1,x_2,Nx_3]+[x_1,Nx_2,Nx_3].
\end{eqnarray}

\begin{defi}
A linear operator $N:\frkg\to \frkg$ is called a Nijenhuis operator if and only if \eqref{eq:Nijenhuis3} and \eqref{eq:Nijenhuis4} hold.
\end{defi}

We have seen that any trivial deformation produces a Nijenhuis operator.
Conversely, any Nijenhuis operator gives a trivial deformation as follows.

\begin{thm}Let $N$ be a Nijenhuis operator for $\frkg$. Then a deformation of $\g$ can be obtained by putting
$$\omega(x_1,x_2,x_3)=[Nx_1,x_2,x_3]+[x_1,Nx_2,x_3]+[x_1,x_2,Nx_3]-N[x_1,x_2,x_3].$$
Furthermore, this deformation is a trivial one.
\end{thm}
\pf It is clear that $\omega=dN$ and $d\omega=ddN=0$. Thus
$\omega$ is a 3-cocycle of $\frkg$ with coefficients in the adjoint representation.
Now we check the fundamental identity \eqref{eq:Lts03} holds for $\omega$. Denote by
\begin{eqnarray*}
J(x_1,x_2,y_1,y_2,y_3)&=&[x_1, x_2, [y_1, y_2, y_3]]- [[ x_1, x_2, y_1], y_2, y_3]\\
&&- [y_1, [ x_1, x_2, y_2], y_3] - [y_1, y_2, [ x_1, x_2, y_3]],\\
J^\omega(x_1,x_2,y_1,y_2,y_3)&=&\omega(x_1, x_2, \omega(y_1, y_2, y_3))- \omega(\omega( x_1, x_2, y_1), y_2, y_3)\\
&&- \omega(y_1, \omega( x_1, x_2, y_2), y_3) - \omega(y_1, y_2, \omega( x_1, x_2, y_3)).
\end{eqnarray*}
A direct computation shows that
\begin{eqnarray*}
&&J^{\omega}(x_1,x_2,y_1,y_2,y_3)\\
&=&J(Nx_1,Nx_2,y_1,y_2,y_3)+N^2J(x_1,x_2,y_1,y_2,y_3)\\
&&-[x_1,x_2,[Ny_1,Ny_2,y_3]+[Ny_1,y_2,Ny_3]+[y_1,Ny_2,Ny_3]-N\omega(y_1,y_2,y_3)]\\
&&+[[Nx_1,Nx_2,y_1]+[Nx_1,x_2,Ny_1]+[x_1,Nx_2,Ny_1]-N\omega(x_1,x_2,y_1),y_2,y_3]\\
&&+[y_1,[Nx_1,Nx_2,y_2]+[Nx_1,x_2,Ny_2]+[x_1,Nx_2,Ny_2]-N\omega(x_1,x_2,y_2),y_3]\\
&&+[y_1,y_2,[Nx_1,Nx_2,y_3]+[Nx_1,x_2,Ny_3]+[x_1,Nx_2,Ny_3]-N\omega(x_1,x_2,y_3)].
\end{eqnarray*}
Therefore we have $J^{\omega}=0$ by the fundamental identity of $\g$ and the Nijenhuis operator condition.
\qed

Note that in the proof of above Theorem we have not used condition \eqref{eq:Nijenhuis3}.
But this condition is important to us since only in this case the $k$'s power of a Nijenhuis operator is also a Nijenhuis operator.

\begin{lem}\label{lem:Nk1}
Let $N$ be a Nijenhuis operator. Then for any $k>0$, we have
\begin{eqnarray}
\label{eq:N1}[x_1,x_2,x_3]_{N^{k+1}}=([x_1,x_2,x_3]_{N^{k}})_{N}.
\end{eqnarray}
\end{lem}

By the above Lemma, we get
\begin{eqnarray*}
([x_1,x_2,x_3]_{N^k})_{N^r}&=&(([x_1,x_2,x_3]_{N^k})_N)_{N^{r-1}}=([x_1,x_2,x_3]_{N^{k+1}})_{N^{r-1}}\\
&=&([x_1,x_2,x_3]_{N^{k+2}})_{N^{r-2}}=\cdots=[x_1,x_2,x_3]_{N^{k+r}}.
\end{eqnarray*}
\begin{lem}\label{lem:Nkr}
Let $N$ be a Nijenhuis operator. Then for any $k,r>0$, we have
\begin{eqnarray}
\label{eq:Nkr}[x_1,x_2,x_3]_{N^{k+r}}=([x_1,x_2,x_3]_{N^{k}})_{N^{r}}.
\end{eqnarray}
\end{lem}

\begin{pro}\label{prop:Nk}
Let $N$ be a Nijenhuis operator. Then for any $k>0$, $N^k$  is also a Nijenhuis operator.
\end{pro}

\pf We prove by introduction. The Proposition is valid for $k=1$. Assume
$$N^k[x_1,x_2,x_3]_{N^k}=[N^kx_1,N^kx_2,x_3]+[N^kx_1,x_2,N^kx_3]+[x_1,N^kx_2,N^kx_3],$$
then we have
\begin{eqnarray*}
&&N^{k+1}[x_1,x_2,x_3]_{N^{k+1}}\\
&=&N^{k}N(([x_1,x_2,x_3]_{N^k})_N)\\
&=&N^k([Nx_1,Nx_2,x_3]_{N^k}+[Nx_1,x_2,Nx_3]_{N^k}+[x_1,Nx_2,Nx_3]_{N^k})\\
&=&[N^{k+1}x_1,N^{k+1}x_2,x_3]+[N^{k+1}x_1,Nx_2,N^kx_3]+[Nx_1,N^{k+1}x_2,N^kx_3]\\
&&+[N^{k+1}x_1,N^kx_2,Nx_3]+[N^{k+1}x_1,x_2,N^{k+1}x_3]+[Nx_1,N^kx_2,N^{k+1}x_3]\\
&&+[N^kx_1,N^{k+1}x_2,Nx_3]+[N^kx_1,Nx_2,N^{k+1}x_3]+[x_1,N^{k+1}x_2,N^{k+1}x_3]\\
&=&[N^{k+1}x_1,N^{k+1}x_2,x_3]+[N^{k+1}x_1,x_2,N^{k+1}x_3]+[x_1,N^{k+1}x_2,N^{k+1}x_3]\\
&&+[N^{k+1}x_1,Nx_2,N^kx_3]+[Nx_1,N^{k+1}x_2,N^kx_3]+[N^{k+1}x_1,N^kx_2,Nx_3]\\
&&+[Nx_1,N^kx_2,N^{k+1}x_3]+[N^kx_1,N^{k+1}x_2,Nx_3]+[N^kx_1,Nx_2,N^{k+1}x_3].
\end{eqnarray*}
The terms in the last two line are zero by \eqref{eq:Nijenhuis3}. Thus Proposition \ref{prop:Nk} is valid for $k+1$.
\qed

Two Nijenhuis operator $N_1$ and $N_2$ are said to be compatible if $N_1+N_2$ is also a Nijenhuis operator.

\begin{pro}\label{prop:N12}
Let $N_1$ and $N_2$ be two Nijenhuis operators. Then they are compatible if and only if
\begin{eqnarray}
 \nonumber  && N_1[x_1,x_2,x_3]_{N_2}+N_2[x_1,x_2,x_3]_{N_1}\\
 \nonumber &=& [N_2x_1,N_1x_2,x_3]+[N_2x_1,x_2,N_1x_3]+[x_1,N_2x_2,N_1x_3]\\
\label{eq:N12-1}&&+[N_1x_1,N_2x_2,x_3]+[N_1x_1,x_2,N_2x_3]+[x_1,N_1x_2,N_2x_3],
\end{eqnarray}
and
\begin{eqnarray}
 \nonumber   &&[N_1x_1,N_1x_2,N_2x_3]+[N_1x_1,N_2x_2,N_2x_3]\\
\label{eq:N12-2} &&+[N_2x_1,N_1x_2,N_1x_3]+[N_2x_1,N_2x_2,N_1x_3]=0.
\end{eqnarray}
\end{pro}

\emptycomment{
\pf By definition, $N_1+N_2$ is a Nijenhuis operator if and only if
\begin{eqnarray*}
&&(N_1+N_2)([x_1,x_2,x_3]_{N_1+N_2})\\
&=&[(N_1+N_2)x_1,(N_1+N_2)x_2,x_3]+[(N_1+N_2)x_1,x_2,(N_1+N_2)x_3]\\
                 &&+[x_1,(N_1+N_2)x_2,(N_1+N_2)x_3],
\end{eqnarray*}
and
\begin{eqnarray*}
[(N_1+N_2)x_1,(N_1+N_2)x_2,(N_1+N_2)x_3]=0.
\end{eqnarray*}
Now it is easy to see the above condition is equivalent to \eqref{eq:N12-1} and \eqref{eq:N12-2}.
\qed
}

\begin{lem}\label{lem:Njk}
Let $N$ be a Nijenhuis operator. Then for any $j,k>0$, we have
\begin{eqnarray}
 \nonumber && N^j[x_1,x_2,x_3]_{N^k}+N^k[x_1,x_2,x_3]_{N^j}\\
 \nonumber &=& [N^kx_1,N^jx_2,x_3]+[N^kx_1,x_2,N^jx_3]+[x_1,N^kx_2,N^jx_3]\\
\label{eq:Njk}&&+[N^jx_1,N^kx_2,x_3]+[N^jx_1,x_2,N^kx_3]+[x_1,N^jx_2,N^kx_3].
\end{eqnarray}
\end{lem}

\pf If $j>k$, then by Proposition \ref{prop:Nk} and Lemma \ref{lem:Nkr} we have
\begin{eqnarray*}
&&N^j[x_1,x_2,x_3]_{N^k}+N^k[x_1,x_2,x_3]_{N^j}\\
&=&N^{j-k}(N^k[x_1,x_2,x_3]_{N^k})+N^k(([x_1,x_2,x_3]_{N^{j-k}})_{N^k})\\
&=&N^{j-k}([N^kx_1,N^kx_2,x_3]+[N^kx_1,x_2,N^kx_3]+[x_1,N^kx_2,N^kx_3])\\
&&+[N^kx_1,N^kx_2,x_3]_{N^{j-k}}+[N^kx_1,x_2,N^kx_3]_{N^{j-k}}+[x_1,N^kx_2,N^kx_3]_{N^{j-k}}\\
&=&N^{j-k}([N^kx_1,N^kx_2,x_3]+[N^kx_1,x_2,N^kx_3]+[x_1,N^kx_2,N^kx_3])\\
&&+[N^jx_1,N^kx_2,x_3]+[N^kx_1,N^jx_2,x_3]+[N^kx_1,N^kx_2,N^{j-k}x_3]-N^{j-k}[N^kx_1,N^kx_2,x_3]\\
&&+[N^jx_1,x_2,N^kx_3]+[N^kx_1,N^{j-k}x_2,N^kx_3]+[N^kx_1,x_2,N^{j}x_3]-N^{j-k}[N^kx_1,x_2,N^kx_3]\\
&&+[N^{j-k}x_1,N^kx_2,N^kx_3]+[x_1,N^jx_2,N^kx_3]+[x_1,N^kx_2,N^{j}x_3]-N^{j-k}[x_1,N^kx_2,N^kx_3]\\
&=& [N^kx_1,N^jx_2,x_3]+[N^kx_1,x_2,N^jx_3]+[x_1,N^kx_2,N^jx_3]\\
&&+[N^jx_1,N^kx_2,x_3]+[N^jx_1,x_2,N^kx_3]+[x_1,N^jx_2,N^kx_3].
\end{eqnarray*}
The case of $j<k$ can be proved similarly and the case of $j=k$ is by Proposition \ref{prop:Nk}.
\qed

Let $N_1=N^j$ and $N_2=N^k$. Then condition \eqref{eq:N12-2} in Proposition \ref{prop:N12} is satisfied by \eqref{eq:Nijenhuis3}
and condition \eqref{eq:N12-1} is satisfied by Lemma \ref{lem:Njk}, thus we get
\begin{pro}
Let $N$ be a Nijenhuis operator. Then for any $j,k>0$, $N^j$ and $N^k$ are compatible.
\end{pro}

It is easy to see that if $N$ is a Nijenhuis operator, then $cN$ is also a Nijenhuis operator, where $c$ is any constant.
Now by Proposition \ref{prop:N12} and Lemma \ref{lem:Njk} we have


\begin{thm}
Let $N$ be a Nijenhuis operator. Then for any polynomial
$P(X)=\sum_{i=1}^n c_i X^i$, the operator $P(N)$ is also a Nijenhuis operator.
\end{thm}

\section{Abelian Extensions of Lie triple systems}
\label{sec:4}

In this section, we study abelian extensions of Lie triple systems.
We show that associated to any abelian extension, there is a representation and a 3-cocycle.
Furthermore, abelian extensions can be classified by the third cohomology group.

An ideal of a Lie triple system $T$ is a subspace $I$ such that $[I,T,T]\subseteq I$.
An ideal $I$ of a Lie triple system $T$ is called an abelian ideal if moreover $[T,I,I]=0$.
Notice that $[T,I,I]=0$ implies that $[I,T,I]=0$ and $[I,I,T]=0$.

\begin{defi}
 Let $(\frkg, [\cdot,\cdot,\cdot])$, $(\frkh, [\cdot,\cdot,\cdot]_\frkh)$,
 $(\hat{\frkg}, [\cdot,\cdot,\cdot]_{\hat{\g}})$ be Lie triple systems and
$i:\frkh\to\hat{\frkg},~~p:\hat{\frkg}\to\frkg$
be homomorphisms. The following sequence of Lie triple systems is a
short exact sequence if $\mathrm{Im}(i)=\mathrm{Ker}(p)$,
$\mathrm{Ker}(i)=0$ and $\mathrm{Im}(p)=\g$,
\begin{equation}\label{diagram:exact}
 \xymatrix{
   0  \ar[r]^{} & \frkh \ar[r]^{i} & \hat{\frkg} \ar[r]^{p} & \frkg  \ar[r]^{} & 0. \\
  }
\end{equation}
In this case, we call $\hat{\frkg}$  an extension of $\frkg$ by
$\frkh$, and denote it by $\E_{\hat{\g}}$.
It is called an abelian extension if $\frkh$ is an abelian ideal of $\hat{\frkg}$, i.e.
$[u,v,\cdot]_{\hat{\g}}=[u,\cdot,v]_{\hat{\g}}=[\cdot,u,v]_{\hat{\g}}=0$,
for all $u,v\in \frkh$.
\end{defi}

A section $\sigma:\frkg\to\hat{\frkg}$ of $p:\hat{\frkg}\to\frkg$
consists of linear maps
$\sigma:\frkg\to\hat{\g}$
 such that  $p\circ\sigma=\id_{\frkg}$.

\begin{defi}
 Two extensions of Lie triple system
 $\E_{\hat{\g}}:0\to\frkh\stackrel{i}{\to}\hat{\g}\stackrel{p}{\to}\g\to0$
 and $\E_{\tilde{\g}}:0\to\frkh\stackrel{j}{\to}\tilde{\g}\stackrel{q}{\to}\g\to0$ are equivalent,
 if there exists a Lie triple system homomorphism $F:\hat{\g}\to\tilde{\g}$  such that the following diagram commutes
\begin{equation}\label{diagram:equivalent}
\xymatrix{
   0  \ar[r]^{} & \frkh \ar[d]_{\id} \ar[r]^{i} & \hat{\frkg} \ar[d]_{F} \ar[r]^{p} & \frkg \ar[d]_{\id} \ar[r]^{} & 0 \\
   0 \ar[r]^{} & \frkh \ar[r]^{j} & \tilde{\frkg} \ar[r]^{q} & \frkg \ar[r]^{} & 0.
   }
\end{equation}
The set of equivalent classes of extensions of $\g$ by $\h$ is denoted by $\Ext(\frkg,\frkh)$.
\end{defi}

Let $\hat{\frkg}$ be an abelian extension of $\frkg$ by
$\frkh$, and $\sigma:\g\to\hat{\g}$ be a section. Define maps from $\ot^2\frkg$ to $\End(\frkh)$ by
\begin{eqnarray}\label{eq:rep}
&&\DD(x)(u)=\DD(x_1,x_2)(u)\triangleq[\sigma(x_1),\sigma(x_2),u]_{\hat{\g}},\\
&&\theta(x)(u)=\theta(x_1,x_2)(u)\triangleq[u,\sigma(x_1),\sigma(x_2)]_{\hat{\g}}.
\end{eqnarray}
Then by antisymmetricity and the Jacobi identity of $[\cdot,\cdot,\cdot]_{\hat{\g}}$ in $\widehat{T}$ we have
\begin{eqnarray*}
\DD(x_1,x_2)(u)=\theta(x_2,x_1)(u)-\theta(x_1,x_2)(u),
\end{eqnarray*}
for all $x=(x_1,x_2)\in\ot^2\frkg$, $u\in\frkh$.

\begin{lem}\label{pro:2-modules}
With the above notations, $(V,\theta)$ is a representation of $\g$ and does not depend on the choice of the section $\sigma$.
Moreover,  equivalent abelian extensions give the same representation.
\end{lem}

\pf
First, we show that $\theta$ is independent of
the choice of $\sigma$. In fact, if we choose another section $\sigma':\g\to \widehat{T}$, then
$$p(\sigma(x_i)-\sigma'(x_i))=x_i-x_i=0
\Longrightarrow\sigma(x_i)-\sigma'(x_i)\in \h\Longrightarrow\sigma'(x_i)=\sigma(x_i)+u_i$$
for some $u\in\h$.

Since we  have  $[u,v,\cdot]_{\hat{\g}}=0=[u,\cdot,v]_{\hat{\g}}$
for all $u,v\in\frkh$, this implies that
\begin{eqnarray*}
[v,\sigma'(x_1),\sigma'(x_2)]_{\hat{\g}}&=&[v,\sigma(x_1)+u_1,\sigma(x_2)+u_2]_{\hat{\g}}\\
&=&[v,\sigma(x_1),\sigma(x_2)+u_2]_{\hat{\g}}+[v,u_1,\sigma(x_2)+u_2]_{\hat{\g}}\\
&=&[v,\sigma(x_1),\sigma(x_2)]_{\hat{\g}}+[v,u_1,\sigma(x_2)]_{\hat{\g}}\\
&=&[v,\sigma(x_1),\sigma(x_2)]_{\hat{\g}}.
\end{eqnarray*}
Thus $\theta$ is independent on the choice of $\sigma$.

Second, we show that $(V,\theta)$ is a representation of $\g$.

By the equality
\begin{eqnarray*}
&&[\si x_1,\si x_2, [u,\si y_1,\si y_2]_{\hat{\g}}]_{\hat{\g}} \\
&=& [[\si x_1,\si x_2,u]_{\hat{\g}},\si y_1,\si y_2]_{\hat{\g}} + [u, [\si x_1,\si x_2,\si y_1]_{\hat{\g}},\si y_2]_{\hat{\g}}+ [u,\si y_1, [\si x_1,\si x_2,\si y_2]_{\hat{\g}}]_{\hat{\g}},
\end{eqnarray*}
we have
\begin{eqnarray*}
\DD(x_1,x_2)\theta(y_1,y_2)(u)=\theta(y_1,y_2)\DD(x_1,x_2)(u)+\theta((x_1,x_2)\circ (y_1,y_2))(u),
\end{eqnarray*}
where we use the fact that
$$\si[x_1,x_2,y_1]_{\g}-[\si x_1,\si x_2,\si y_1]_{\hat{\g}}\in V\cong \mathrm{Ker}(p),$$
and that $\frkh$ is an abelian ideal of $\hat{\frkg}$,
$$[u,\si[x_1,x_2,y_1]_{\g}-[\si x_1,\si x_2,\si y_1]_{\hat{\g}},\si y_2]_{\hat{\g}}=0.$$
Thus we obtain the condition (R1).

By the equality
\begin{eqnarray*}
&&[u,\si x_1, [\si y_1,\si y_2,\si y_3]_{\hat{\g}}]_{\hat{\g}} \\
&=& [[u,\si x_1,\si y_1]_{\hat{\g}},\si y_2,\si y_3]_{\hat{\g}} + [\si y_1, [u,\si x_1,\si y_2]_{\hat{\g}},\si y_3]_{\hat{\g}}+ [\si y_1,\si y_2, [u,\si x_1,\si y_3]_{\hat{\g}}]_{\hat{\g}},
\end{eqnarray*}
we have
\begin{eqnarray*}
\theta(x_1,[y_1, y_2, y_3])(u)= \theta (y_2, y_3)\theta(x_1,y_1)(u) - \theta (y_1, y_3)\theta(x_1,y_2)(u) + \DD (y_1, y_2)\theta(x_1,y_3)(u).
\end{eqnarray*}
Thus we get the condition (R2). Therefore we deduce that $(V,\theta)$ is a representation of $\g$.

At last, suppose that $\E_{\hat{\g}}$ and $\E_{\tilde{\g}}$ are equivalent abelian extensions, and $F:\hat{\g}\to\tilde{\g}$ is the Lie triple system homomorphism satisfying $F\circ i=j$, $q\circ F=p$.
Choosing linear sections $\sigma$ and $\sigma'$ of $p$ and $q$, we get $qF\sigma(x_i)=p\sigma(x_i)=x_i=q\sigma'(x_i)$,
then $F\sigma(x_i)-\sigma'(x_i)\in \Ker (q)\cong\h$. Thus, we have
$$
[u,\sigma(x_1),\sigma(x_2)]_{\hat{\g}}=[u,F\sigma(x_1),F\sigma(x_2)]_{\tilde{\g}}=[u,\sigma'(x_1),\sigma'(x_2)]_{\tilde{\g}}.
$$
Therefore, equivalent abelian extensions give the same $\theta$. The proof is finished.
\qed\vspace{3mm}

Let $\sigma:\frkg\to\hat{\frkg}$  be a
section of the abelian extension. Define the following map:
\begin{equation}\label{eq:coc}
\omega(x_1,x_2,x_3)\triangleq[\sigma(x_1),\sigma(x_2),\sigma(x_3)]_{\hat{\frkg}}-\sigma([x_1,x_2,x_3]_{\frkg}),
\end{equation}
for all $x_1,x_2,x_3\in\frkg$.

\begin{lem}\label{thm:2-cocylce}
Let $0\to\frkh{\to}\hat{\g}{\to}\g\to 0$ be an abelian extension of $\g$ by $\h$.
Then $\omega$ defined by \eqref{eq:coc} is a 3-cocycle of $\g$ with coefficients in $\frkh$,
where the representation $\theta$ is given by \eqref{eq:rep}.
\end{lem}

\pf  It is easy to see that $\omega$ defined above is antisymmetric in the first two variables and satisfies the Jacobi identity.
By the equality
\begin{eqnarray*}
&&[\si x_1,\si x_2, [\si y_1,\si y_2,\si y_3]_{\hat{\g}}]_{\hat{\g}} \\
&=& [[\si x_1,\si x_2,\si y_1]_{\hat{\g}},\si y_2,\si y_3]_{\hat{\g}} + [\si y_1, [\si x_1,\si x_2,\si y_2]_{\hat{\g}},\si y_3]_{\hat{\g}}\\
&& + [\si y_1,\si y_2, [\si x_1,\si x_2,\si y_3]_{\hat{\g}}]_{\hat{\g}},
\end{eqnarray*}
we get that the left hand side is equal to
\begin{eqnarray*}
&&[\si x_1,\si x_2, [\si y_1,\si y_2,\si y_3]_{\hat{\g}}]_{\hat{\g}}\\
&=&[\si x_1,\si x_2, \omega(y_1,y_2,y_3)+\sigma([y_1,y_2,y_3]_\frkg)]_{\hat{\g}}\\
&=&\DD(x_1, x_2)\omega(y_1,y_2,y_3)+[\si x_1,\si x_2,\sigma([y_1,y_2,y_3]_\frkg)]_{\hat{\g}}\\
&=&\DD(x_1, x_2)\omega(y_1,y_2,y_3)+\omega(x_1,x_2,[y_1,y_2,y_3]_\frkg)+\sigma([x_1,x_2,[y_1,y_2,y_3]_\frkg]_\frkg).
\end{eqnarray*}
Similarily, the right hand side is equal to
\begin{eqnarray*}
&& \theta(y_2, y_3)\omega(x_1,x_2,y_1)+\omega([ x_1, x_2, y_1]_\frkg, y_2, y_3+\sigma([[ x_1, x_2, y_1]_\frkg, y_2, y_3]_\frkg)\\
&&-\theta(y_1, y_3)\omega(x_1,x_2,y_2)+\omega(y_1, [x_1, x_2, y_2]_\frkg, y_3)+\sigma([y_1, [ x_1, x_2, y_2]_\frkg, y_3]_\frkg)\\
&&+\DD(y_1, y_2) \omega(x_1,x_2,y_3)+ \omega(y_1, y_2, [ x_1, x_2, y_3]_\frkg)+\sigma([y_1, y_2, [ x_1, x_2, y_3]_\frkg]_\frkg).
\end{eqnarray*}
Thus we have
\begin{eqnarray*}
&& \omega( x_1, x_2,[y_1, y_2, y_3]_\frkg)+\DD(x_1, x_2)\omega(y_1, y_2, y_3)\\
&=&\omega([x_1, x_2, y_1]_\frkg, y_2, y_3) + \omega(y_1, [x_1, x_2,y_2]_\frkg, y_3) + \omega(y_1, y_2, [x_1, x_2,y_3]_\frkg)\\
&&+\theta(y_2, y_3)\omega(x_1,x_2,y_1)- \theta(y_1, y_3)\omega(x_1,x_2,y_2) + \DD(y_1, y_2)\omega(x_1,x_2,y_3).
\end{eqnarray*}
This is exactly the 3-cocycle condition in Definition \ref{def:1coc}.
\qed

Now we can transfer the Lie triple system structure on $\hat{\g}$ to the Lie triple system structure on $\frkg\oplus\frkh$ using the 3-cocycle given above.
More precisely, we have
\begin{lem}
Let $\frkg$ be a Lie triple system, $(V, \theta)$ be a $\frkg$-module and $\omega: \ot^3 \g\to V$ be a 3-cocycle.
Then $\g\oplus V$ is a Lie triple system under the following multiplication:
\begin{eqnarray*}
&&[x_1 + u_1, x_2 + u_2, x_3 + u_3]_\omega\\
&=&[x_1, x_2, x_3] + \omega(x_1, x_2, x_3)+ \DD(x_1, x_2)(u_3)-\theta(x_1, x_3)(u_2) + \theta( x_2, x_3)(u_1),
\end{eqnarray*}
where $x_1, x_2, x_3 \in \frkg$ and $u_1, u_2, u_3 \in V$. This Lie triple system is denoted by $\g\oplus_\omega\h$.
\end{lem}

\pf It is easy to see that
\begin{eqnarray*}
&&[x_1 + u_1, x_1 + u_1, x_2 + u_2]_\omega\\
&=&[x_1, x_1, x_2] + \omega(x_1, x_1, x_2)+ \DD(x_1, x_1)(u_2)-\theta(x_1, x_2)(u_1) + \theta( x_1, x_2)(u_1)\\
&=&0,
\end{eqnarray*}
and
\begin{eqnarray*}
&&[x_1 + u_1, x_2 + u_2, x_3 + u_3]_\omega+[x_2 + u_2, x_3 + u_3, x_1 + u_1]_\omega+[x_3 + u_3, x_1 + u_1, x_2 + u_2]_\omega\\
&=&[x_1, x_2, x_3] + \omega(x_1, x_2, x_3)+ \DD(x_1, x_2)(u_3)-\theta(x_1, x_3)(u_2) + \theta( x_2, x_3)(u_1)\\
&&+[x_2, x_3, x_1] + \omega(x_2, x_3, x_1)+ \DD(x_2, x_3)(u_1)-\theta(x_2, x_1)(u_3) + \theta( x_3, x_1)(u_2)\\
&&+[x_3, x_1, x_2] + \omega(x_3, x_1, x_2)+ \DD(x_3, x_1)(u_2)-\theta(x_3, x_2)(u_1) + \theta( x_1, x_2)(u_3)\\
&=&0,
\end{eqnarray*}
where the representation terms related to $\DD$ and $\theta$ cancel out and the remaining terms are equal to zero since
$\omega$ is a 3-cocycle.

Now it suffices to verify the fundamental identity.
\begin{eqnarray*}
&&[x_1 + u_1, x_2 + u_2,[y_1 + v_1, y_2 + v_2, y_3 + v_3]_\omega]_\omega\\
&=&[x_1 + u_1, x_2 + u_2, [y_1, y_2, y_3] + \omega(y_1, y_2, y_3)+ \DD(y_1, y_2)(v_3)-\theta(y_1, y_3)(v_2) + \theta( y_2, y_3)(v_1)]\\
&=&[x_1, x_2, [y_1, y_2, y_3]] + \omega(x_1, x_2, [y_1, y_2, y_3])-\theta(x_1, [y_1, y_2, y_3])(u_2) + \theta( x_2, [y_1, y_2, y_3])(u_1)\\
&& +\DD(x_1, x_2)(\omega(y_1, y_2, y_3)+ \DD(y_1, y_2)(v_3)-\theta(y_1, y_3)(v_2) + \theta( y_2, y_3)(v_1)),
\end{eqnarray*}
\begin{eqnarray*}
&&[[x_1 + u_1, x_2 + u_2,y_1 + v_1]_\omega, y_2 + v_2, y_3 + v_3]_\omega\\
&=&[[x_1, x_2, y_1] + \omega(x_1, x_2, y_1)+ \DD(x_1, x_2)(v_1)-\theta(x_1, y_1)(u_2) + \theta(x_2, y_1)(u_1), y_2 + v_2, y_3 + v_3]\\
&=&[[x_1, x_2, y_1], y_2, y_3] + \omega([x_1, x_2, y_1], y_2, y_3)+\DD([x_1, x_2, y_1],y_2)(v_3) - \theta([x_1, x_2, y_1],y_3)(v_2)\\
&& +\theta(y_2, y_3)(\omega(x_1, x_2, y_1)+ \DD(x_1, x_2)(v_1)-\theta(x_1, y_1)(u_2) + \theta(x_2, y_1)(u_1)),
\end{eqnarray*}
\begin{eqnarray*}
&&[y_1 + v_1, [x_1 + u_1, x_2 + u_2, y_2 + v_2]_\omega, y_3 + v_3]_\omega\\
&=&[y_1 + v_1, [x_1, x_2, y_2] + \omega(x_1, x_2, y_2)+ \DD(x_1, x_2)(v_2)-\theta(x_1, y_2)(u_2) + \theta( x_2, y_2)(u_1), y_3 + v_3]\\
&=&[y_1, [x_1, x_2, y_2],y_3] + \omega(y_1, [x_1, x_2, y_2], y_3)+\DD(y_1, [x_1, x_2, y_2])(v_3) + \theta([x_1, x_2, y_2], y_3)(v_1)\\
&& -\theta(y_1, y_3)(\omega(x_1, x_2, y_2)+ \DD(x_1, x_2)(v_2)-\theta(x_1, y_2)(u_2) + \theta( x_2, y_2)(u_1)),
\end{eqnarray*}
\begin{eqnarray*}
&&[y_1 + v_1, y_2 + v_2,[x_1 + u_1, x_2 + u_2, y_3 + v_3]_\omega]_\omega\\
&=&[y_1 + v_1, y_2 + v_2, [x_1, x_2, y_3] + \omega(x_1, x_2, y_3)+ \DD(x_1, x_2)(v_3)-\theta(x_1, y_3)(u_2) + \theta( x_2, y_3)(u_1)]\\
&=&[y_1, y_2, [x_1, x_2, y_3]] + \omega(y_1, y_2, [x_1, x_2, y_3])-\theta(y_1, [x_1, x_2, y_3])(v_2) + \theta( y_2, [x_1, x_2, y_3])(v_1)\\
&& +\DD(y_1, y_2)(\omega(x_1, x_2, y_3)+ \DD(x_1, x_2)(v_3)-\theta(x_1, y_3)(u_2) + \theta( x_2, y_3)(u_1)),
\end{eqnarray*}
It follows that
\begin{eqnarray*}
&&[x_1 + u_1, x_2 + u_2,[y_1 + v_1, y_2 + v_2, y_3 + v_3]_\omega]_\omega\\
&=&[[x_1 + u_1, x_2 + u_2,y_1 + v_1]_\omega, y_2 + v_2, y_3 + v_3]_\omega\\
&&+[y_1 + v_1, [x_1 + u_1, x_2 + u_2, y_2 + v_2]_\omega, y_3 + v_3]_\omega\\
&&+[y_1 + v_1, y_2 + v_2,[x_1 + u_1, x_2 + u_2, y_3 + v_3]_\omega]_\omega
\end{eqnarray*}
by the fact that $(V, \theta)$ is a represention and by the 3-cocycle condition \eqref{eq:1coc03}.
\qed

\begin{lem}\label{thm:2-cocylce}
 Two abelian extensions of Lie triple systems $0\to\frkh{\to}\g\oplus_\omega\h{\to}\g\to 0$
 and  $0\to\frkh{\to}\g\oplus_{\omega'}\h{\to}\g\to 0$ are equivalent if and only if $\omega$ and $\omega'$ are in the same cohomology class.
\end{lem}

\pf Let $F:\g\oplus_\omega\h\to \g\oplus_{\omega'}\h$ be the corresponding homomorphism. Then we have
\begin{eqnarray}\label{hhh}
&&F[x_1, x_2, x_3]_{\omega}=[F(x_1),F(x_2),F(x_3)]_{\omega'}.
\end{eqnarray}
Since $F$ is an equivalence of extensions, there exists $\nu:\frkg\to \frkh$ such that
\begin{eqnarray}\label{equiv}
F(x_i+u_i)=x_i+\nu(x_i)+u_i,\quad i=1,2,3.
\end{eqnarray}
The left hand side of \eqref{hhh} is equal to
\begin{eqnarray*}
&&F_1([x_1, x_2, x_3]+\omega(x_1, x_2, x_3))\\
&=&[x_1, x_2, x_3]+\omega(x_1, x_2, x_3)+\nu([x_1, x_2, x_3]),
\end{eqnarray*}
and the right hand side of \eqref{hhh} is equal to
\begin{eqnarray*}
&&[x_1+\nu(x_1),x_2+\nu(x_2),x_3+\nu(x_3)]_{\omega'}\\
&=&[x_1,x_2,x_3]+\omega'(x_1,x_2,x_3)\\
&&+\DD(x_1, x_2)\nu(x_3)-\theta(x_1,x_3)\nu(x_2)+\theta(x_2,x_3)\nu(x_1).
\end{eqnarray*}
Thus we have
\begin{eqnarray}\label{eq:exact4}
\nonumber(\omega-\omega')(x_1, x_2, x_3)&=&\DD(x_1, x_2)\nu(x_3)-\theta(x_1,x_3)\nu(x_2)+\theta(x_2,x_3)\nu(x_1)\\
&&-\nu([x_1, x_2, x_3]),
\end{eqnarray}
that is $\omega-\omega'=d\nu$. Therefore $\omega$ and $\omega'$ are in the same cohomology class.
Conversely, if $\omega$ and $\omega'$  are in the same cohomology class, assume that $\omega-\omega'=d\nu$.
Then we can define $F$ by \eqref{equiv}. Similar as the above proof, we can show that $F$ is an equivalence.
We omit the details.
\qed

\begin{thm}
Let $\frkg$ be a Lie triple system and $(V, \theta)$ be a $\frkg$-module.
Then there is a one-to-one correspondence between equivalence classes of abelian extensions of Lie triple system $\frkg$ by $\h$ and the cohomology group $\mathbf{H}^3(\frkg, V)$.
\end{thm}

For the cohomology group $\mathbf{H}^5(\frkg, V)$, it will be correspondence to equivalence classes of crossed module extensions of Lie triple systems which will defined in another paper.

\section*{Acknowledgments}
The author would like to thank the referee for careful reading of the manuscript and for valuable comments and suggestions
which were very helpful to improve this paper.
Part of this work was done while the author was visiting Courant Research Centre, Georg-August Universit\"{a}t G\"{o}ttingen from June to September 2013.
He is grateful to the University of G\"{o}ttingen for hospitality and financial support.


\vskip7pt

\footnotesize{\noindent
 College of Mathematics and Information Science\\
 Henan Normal University\\
 Xinxiang 453007, P. R. China\\
 zhangtaozata@gmail.com}


\begin{thebibliography}{99}


\bibitem{Ca}
Cartan, E.,  "Oeuvres completes", Part 1, vol. 2, nos. 101, 138, Gauthier-Villars, Paris, 1952.


\bibitem{CLP}
Casas, J. M., J.-L. Loday and T. Pirashvili, \emph{Leibniz $n$-algebras}, Forum Math. {\bf 14} (2002), 189--207.

\bibitem{DT}
Daletskii, Y. I.,and L. A. Takhtajan, \emph{Lie and Lie algebra structures for Nambu algebra},
Lett. Math. Phys. {\bf39} (1997), 127--141.


\bibitem{D} Dorfman, I., ''Dirac Structures and Integrability of Nonlinear Evolution Equation'',
John Wiley \& Sons, Ltd., Chichester, 1993.


\bibitem{Gau} Gautheron, P., \emph{Some remarks concerning Nambu mechanics}, Lett. Math. Phys. {\bf 37} (1996), 103--116.




\bibitem{Ha}
Harris, B., \emph{Cohomology of Lie triple systems and Lie algebras with involution}, Trans. Amer. Math. Soc. {\bf98} (1961), 148--162.

\bibitem{HP}
Hodge, T. L., and B. J. Parshall,
\emph{On the representation theory of Lie triple systems}, Trans. Amer. Math. Soc. {\bf354} (2002), 4359--4391.

\bibitem{Ja1}
Jacobson, N., \emph{Lie and Jordan triple systems}, Amer. J. Math., {\bf71} (1949), 149--170.



\bibitem{KM} Kosmann-Schwarzbach, Y., and F. Magri, \emph{Poisson-Nijenhuis structures},
Ann. Inst. H. Poincar\'e Phys. Th\'eor.  {\bf53} (1990), 35--81.



\bibitem{KT}
Kubo, F., and Y. Taniguchi, \emph{A controlling cohomology of the deformation theory of Lie triple systems}, J. Algebra
{\bf278} (2004), pp. 242--250.


\bibitem{LWD}
Lin, J.,  Y. Wang and S. Deng, \emph{T*-extension of Lie triple systems}, Linear Algebra Appl. {\bf431} (2009), 2071--2083.


\bibitem{Li}
Lister, W.G., \emph{A structure theory of Lie triple systems}, Trans. Amer. Math. Soc. {\bf72} (1952) 217--242.




\bibitem{Loday}
Loday, J.-L., and T. Pirashvili, \emph{Universal enveloping algebras of Leibniz algebras and (co)-homology},
Math. Ann. {\bf296} (1993), 139--158.



\bibitem{NR}
Nijenhuis, A., and R. W. Richardson, \emph{Cohomology and deformations in graded Lie algebras},
Bull. Amer. Math. Soc. {\bf72} (1966), 1--29.


\bibitem{Ya}
Yamaguti, K., \emph{On the cohomology space of Lie triple systems}, Kumamoto J. Sci. A. {\bf5} (1960) 44--52.
\end{thebibliography}
\end{document}